\documentclass[12pt]{article}

\usepackage{mathtools}
\usepackage{microtype}
\usepackage[usenames,dvipsnames,svgnames,table]{xcolor}
\usepackage{amsfonts,amsmath, amssymb,amsthm,amscd}
\usepackage[height=9in,width=6.5in,margin=1in]{geometry}
\usepackage{verbatim}
\usepackage{tikz}
\usepackage{latexsym}
\usepackage[pdfauthor={Lyons, Peres, and Sun},bookmarksnumbered,hyperfigures,colorlinks=true,citecolor=BrickRed,linkcolor=BrickRed,urlcolor=BrickRed,pdfstartview=FitH]{hyperref}
\usepackage{mathrsfs}
\usepackage{graphicx}

\newcommand{\old}[1]{}
\newtheorem{theorem}{Theorem}[section]
\newtheorem{corollary}[theorem]{Corollary}
\newtheorem{lemma}[theorem]{Lemma}
\newtheorem{proposition}[theorem]{Proposition}

\newtheorem{definition}[theorem]{Definition}

\numberwithin{equation}{section}

\newcommand{\be}{\begin{equation}}
	\newcommand{\ee}{\end{equation}}
\def\ba{\begin{align}}
	\def\ea{\end{align}}

\newcommand{\eff}{\mathrm{R}_{\mathrm{eff}}}
\newcommand{\effw}{\mathrm{R}_{\textrm{w-eff}}}

\newcommand{\tree}{\mathcal T_o}

\newcommand{\ray}{\mathrm{Ray}}
\newcommand{\bush}{\mathrm{Bush}}
\renewcommand{\path}{\mathcal{P}}
\newcommand{\LE}{\mathbf{LE}}

\newcommand{\graph}{\overline{\mathcal{T}_o}}

\newcommand{\N}{\mathbb{N}}

\newcommand{\Z}{\mathbb{Z}}
\newcommand{\R}{\mathbb{R}}

\newcommand\st{\,;\;}

\renewcommand{\P}{\mathbb{P}}
\newcommand{\E}{\mathbb{E}}

\newcommand\xor{\mathbin{\triangle}}
\newcommand{\notion}[1]{{\bf  \textit{#1}}}
\newcommand{\op}[1]{\operatorname{#1}}

\DeclareMathOperator{\WSF}{\mathsf{WSF}}
\DeclareMathOperator{\UST}{\mathsf{UST}}
\DeclareMathOperator{\FSF}{\mathsf{FSF}}
\DeclareMathOperator{\USF}{\mathsf{USF}}

\def\WSFs{\mathfrak{F}_{\mathrm w}}
\def\FSFs{\mathfrak{F}_{\mathrm f}}
\def\USTs{\mathfrak{T}}
\def\iw{i_{\mathrm w}}

\def\cT{\mathcal{T}}

\def\cP{\mathcal{P}}

\def\cN{\mathcal{N}}

\def\cI{\mathcal{I}}

\def\cG{\mathcal{G}}

\def\cE{\mathcal{E}}

\def\cC{\mathcal{C}}

\def\cA{\mathcal{A}}

\newcommand{\defeq}{:=}
\newcommand{\ol}{\overline}
\newcommand{\wt}{\widetilde}

\newcommand{\wh}{\widehat}

\newcommand\Seq[1]{\langle #1 \rangle}

\usepackage[normalem]{ulem}

\providecommand{\keywords}[1]
{
	\small	
	\textbf{\textit{Keywords---}} #1
}

\begin {document}
\title{Induced graphs of uniform spanning  forests}
\author{
	Russell Lyons\thanks{Department of Mathematics, Indiana
		University. Partially supported by the National
		Science Foundation under grants DMS-1007244 and DMS-1612363 and by Microsoft Research.
		Email: \protect\url{rdlyons@indiana.edu}.}
	\and
	Yuval Peres\thanks{Kent State University, Kent, OH.
		Email: \protect\url{yuval@yuvalperes.com}.}
	\and
	Xin Sun\thanks{Department of Mathematics, Columbia University. Partially supported by Simons Society of Fellows, by NSF Award DMS-1811092, and by
		Microsoft Research.
		Email: \protect\url{xinsun@math.columbia.edu}.}
}
\date{}
\maketitle
\begin{abstract}
Given a subgraph $H$ of a graph $G$,    the  induced graph of $H$ is the largest subgraph of $G$ whose vertex set is the same as that of $H$. Our paper concerns the induced graphs of the 
components of $\WSF(G)$, the wired uniform spanning forest on $G$, and, to a lesser extent, $\FSF(G)$, the free uniform spanning forest.
We show that the induced graph of each component of  $\WSF(\mathbb Z^d$) is almost surely recurrent when $d\ge 8$. 
Moreover, the effective resistance between two points on the ray of the tree to infinity within a component grows linearly when $d\ge9$. 
For any vertex-transitive graph $G$,  we establish the following resampling property: 
Given a vertex $o$ in $G$, let $\tree$ be the component of $\WSF(G)$ containing $o$ and $\graph$ be its induced graph.   
Conditioned on $\graph$, the tree $\tree$ is distributed  as $\WSF(\graph)$.
For any graph $G$,  we also show that if
$\tree$ is the component of $\FSF(G)$ containing $o$ and $\graph$ is its induced graph,
then conditioned on $\graph$, the tree $\tree$ is distributed  as $\FSF(\graph)$.
\end{abstract}

\begin{abstract}
\'Etant donn\'e un sous-graphe $H$ d'un graphe $G$, le graphe induit de $H$
est le plus grand sous-graphe de $G$ dont l'ensemble de sommets est le
m\^eme que celui de $H$. Notre article concerne les graphes induits des
composants connexes de $\WSF (G)$, la for\^et recouvrante uniforme c\^abl\'ee
sur $G$, et, dans une moindre mesure, $\FSF (G)$, la for\^et recouvrante
uniforme libre.  Nous montrons que le graphe induit de chaque composant de
$\WSF (\mathbb Z^d$) est presque s\^urement r\'ecurrent lorsque $d \ge
8$.  De plus, la r\'esistance effective entre deux points du rayon de
l'arbre \`a l'infini au sein d'un composant cro\^\i t lin\'eairement lorsque
$d \ge 9$.  Pour tout graphe transitif \`a sommets $G$, nous \'etablissons la
propri\'et\'e de r\'e\'echantillonnage suivante: \'Etant donn\'e un sommet
$o$ dans $G$, soit $\tree$ le composant de $\WSF (G)$ qui contient $o$ et
$\graph$ son graphe induit.  Conditionn\'e sur $\graph$, l'arbre $\tree$
est distribu\'e comme $\WSF (\graph)$.  Pour tout graphe $G$, nous
montrons \'egalement que si $\tree$ est le composant de $\FSF (G)$ qui
contient $o$ et $\graph$ est son graphe induit, alors conditionn\'e sur
	$\graph$, l'arbre $\tree$ est distribu\'e comme $\FSF (\graph)$.
\end{abstract}


\section{Introduction}\label{sec:intro}
Given a finite, connected graph $G$, the \notion{uniform spanning tree} (UST)  on $G$, which we denote by $\UST(G)$, is the  uniform measure on the set of spanning trees of $G$.  Given a (locally finite) infinite, connected graph $G$, notions of  ``uniform spanning tree'' can be defined via limiting procedures. Suppose $\Seq{G_n}$ is a sequence of finite, connected subgraphs  of $G$.
We call $\Seq{G_n}$  an  \notion{exhaustion} of $G$ if  $G_n\subset G_{n+1}$ and
$\bigcup G_n= G$.  According to \cite{Pemantle},  given an exhaustion $\Seq{G_n}$ of
$G$ and a fixed, finite subgraph $H$ of $G$, the weak limit of $\UST(G_n)\cap H$
exists. Varying $H$, one obtains a probability measure on subgraphs of $G$, which is called
the \notion{free spanning forest} (FSF) of $G$ and denoted by $\FSF(G)$. On
the other hand, for any finite, connected subgraph $H$ of $G$, let $\wh H$
be the graph obtained by identifying all  vertices not in $H$ to a
single vertex. We call $\UST(\wh H)$  the \notion{wired spanning forest} (WSF)
of $H$ (relative to $G$), which we denote   by $\WSF(H)$. Given an exhaustion $\Seq{G_n}$ of
$G$ consisting of induced subgraphs  and a  fixed, finite subgraph $H$ of $G$,  the weak limit of $\WSF(G_n)\cap H$
exists. Varying $H$, one obtains a probability measure on subgraphs of $G$, which is called
the wired spanning forest of $G$ and denoted by $\WSF(G)$. 

 Neither $\WSF$ and $\FSF$  can have no cycles, but both can have more than one
(connected) component. This  justifies the notion of spanning forest.
$\UST$ and its infinite-volume extensions have been an important object in
probability and mathematical physics for the last three decades. See
\cite{USF,TreeBook} for a basic reference for some of the theory.
See \cite{Schramm00,LSW-UST,Dubedat-GFF,Shef-burger,Bernoulli-UST} for the important role that UST played in the recent development in two dimensional statistical physics and Schramm--Loewner evolution. 

For infinite graphs, $\WSF$ is much better understood than $\FSF$. In \cite{Wilson}, David Wilson provided an efficient  algorithm to sample $\UST$ on finite graphs. It was soon extended to sample $\WSF$ on infinite graphs \cite{USF} (see also Section~\ref{sec:8D}). This powerful tool allows one to study $\WSF$ directly via simple random walk. In particular, it is proved in \cite{USF} that $\WSF(G)$ is concentrated on the set of forests with a unique component if and only if two simple random walks on $G$ intersect a.s. In contrast, there is no known simple condition to determine whether an $\FSF$ has a unique component. 
If $G$ is a Cayley graph of the $d$-dimensional integer lattice $\Z^d$ for $d\in \N$,  then $\WSF(\Z^d)$ and $\FSF(\Z^d)$ coincide.
Moreover, it was shown in \cite{Pemantle} that $\WSF(\Z^d)$ has a unique component when $1\le d\le 4$ and infinitely many components when $d\ge 5$. 
 In the latter case, the collections of components exhibit an intriguing geometry \cite{GeometryofUSF,Component}. 

If $H$ is a subgraph of a graph  $G$, including the case that $H$ may have no edges, then the \notion{induced subgraph} determined by $H$ is the largest subgraph of $G$ whose vertex set is the same as that of $H$.
If $H$ is a subgraph of a graph  $G$,   we define the \notion{induced-component graph} $\ol H$ of $H$ to be the largest subgraph of $G$ whose vertex set  is the same as that
of $H$ and that has the same connected components as $H$.  Thus, an edge of $G$ belongs to $\ol H$ if and only if both its endpoints belong to
the same component of $H$.
We also have that $\ol H$ is the union of the induced subgraphs determined by the components of $H$.

Before stating our main results, we make the following conventions throughout the paper. 
We will  use $\WSF,\FSF,\UST$ to denote either probability measures or their samples as long as it is clear from the context what we are referring to.  When there is a risk of ambiguity,  we use $\UST(G),\FSF(G),\WSF(G)$ to represent probability measures and   $\USTs(G),\FSFs(G),\WSFs(G)$ to represent their corresponding samples.
Similarly, we write $\ol{\WSF(G)}$ for either the law of $\ol{\WSFs(G)}$ or for its sample, $\ol{\WSFs(G)}$, and likewise for the free versions.

The main object of interest in this paper  is $\ol{\WSF(\Z^d)}$, which reflects the geometry of $\WSF(\Z^d)$ as a subgraph embedded in $\Z^d$.
Since $\ol{\WSF(\Z^d)}=\Z^d$ when $1\le d\le 4$, the only interesting case is when $d\ge 5$.
On the one hand,  components of $\ol{\WSF(\Z^d)}$ have stochastic dimension $4$ for all $d\ge 5$ \cite{GeometryofUSF}.
On the other hand, Morris \cite{Morris}  proved that for every graph $G$,  simple random walk on each component of $\WSF(G)$ is a.s.\ recurrent.
This leads to the intriguing question of  whether the components of $\ol{\WSF(\Z^d)}$ are recurrent or transient.

\begin{theorem}\label{thm:recurrence}
If $d\ge 8$, almost surely each connected component of\/ $\ol{\WSF(\Z^d)}$ is recurrent.
\end{theorem}
For a graph $G=(V,E)$, let $f$ be a real function from $V$ to $\R$, and let 
\begin{equation}\label{eq:energy}
\cE(f) :=\frac12\sum_{x,y\in G \st x\sim y} \bigl(f(x)-f(y) \bigr)^2,
\end{equation}where $x\sim y$ means $x$ and $y$ are adjacent in $G$. 
Given two disjoint subsets
$A$ and $B$ of $V$, the \notion{effective resistance}
between $A$ and $B$ is defined by
\[
\eff^G(A,B):=\Bigl(\inf\bigl\{\cE(f) \st f|_A=1, f|_B=0\bigr\}\Bigr)^{-1}.
\]

\begin{definition}\label{def:rays}
An \notion{end} of a tree is an equivalence classes of infinite simple paths in the tree, where
two paths are equivalent if their symmetric difference is finite.
Given a vertex $v$ in a forest, write $T_v$ for the component of the forest that contains $v$. If $T_v$ has one end, then write $\ray_v = \Seq{\ray_v(n)}_{n\ge 0}$ for the unique infinite, simple path in $T_v$ starting at $v$.
Given $v\in \Z^d$, let $\cT_v$ be the component of $\WSF(\Z^d)$ containing $v$. 
\end{definition}

As discussed following Corollary~\ref{cor:wired2}, every tree $\cT_v$ has one end a.s.

\begin{theorem}\label{thm:linear}
Let $v\in \Z^d$.
If $d\ge 9$, then \(\liminf_{n\to \infty}n^{-1}\eff^{\ol{\cT_v}}\bigl(v,\ray_v(n)\bigr) >0 \) almost surely.
\end{theorem}
Since $\eff^{\ol{\cT_v}}\bigl(v,\ray_v(n)\bigr)\le n$, Theorem~\ref{thm:linear} means $\eff^{\ol{\cT_v}}\bigl(v,\ray_v(n)\bigr)$ grows linearly.

Theorems~\ref{thm:recurrence} and~\ref{thm:linear} leave open the natural
questions whether components of   $\ol{\WSF(\Z^d)}$
are recurrent or transient for $d=5,6,7$ and what the growth rate of $\eff^{\ol{\cT}_v}\bigl(v,\ray_v(n)\bigr)$ is for $d=5,6,7,8$.
Although we do not address those problems here, we prove the following resampling property of $\WSF(\Z^d)$ for all dimensions,
which has implications for the behavior of random walks on the components of $\ol{\WSF(\Z^d)}$.

We will use the following notion. Given a graph $G$, let  $H$ be a random subgraph of $G$ whose components are infinite graphs. 
We write $\WSF (H)$ for the unconditional law of the random subgraph of $G$  obtained by first sampling $H$ and then sampling a  WSF independently on each component of this instance of $H$.  We similarly define $\FSF(H)$.
\begin{theorem}\label{thm:resample-Z}
For  all $d\in\N$, $\>\WSF(\ol{\WSFs(\Z^d)})=\WSF(\Z^d)$, that is, the two measures agree.
\end{theorem}
Theorem~\ref{thm:resample-Z} implies that for each $v\in\Z^d$,  $\>\WSF(\ol{\cT_v})$ a.s.\ has a single component for $\cT_v$ as in Definition~\ref{def:rays}.
Therefore, two independent simple random walks on $\ol{\cT_v}$ a.s.\ intersect.

Theorems~\ref{thm:recurrence} and~\ref{thm:linear} are proved in Section~\ref{sec:8D} and~\ref{sec:highD} respectively, using quantitative arguments. 
A \notion{vertex-transitive graph} is a  graph such that given any two
vertices, there exists a graph automorphism mapping one vertex to the
other (see Section~\ref{subsec:transitive}). 
Theorem~\ref{thm:recurrence} can be extended to all  vertex-transitive graphs whose volume growth is at least $r \mapsto r^8$, while
the argument for Theorem~\ref{thm:linear} works for unimodular vertex-transitive graphs whose volume growth is at least  $r \mapsto r^9$ (see the definition of unimodular preceding Theorem~\ref{thm:linear1}).

On the other hand, Theorem~\ref{thm:resample-Z} is a corollary of the following  set of general results, which will be proved in Section~\ref{sec:resample} by qualitative arguments. 
\begin{theorem}\label{thm:resample}
	For any locally finite,  infinite, connected graph $G$, we have  
	\begin{equation}\label{eq:resample}
		\FSF(\ol{\WSFs(G)})=\WSF(G)\quad \textrm{and}\quad	\FSF(\ol{\FSFs(G)})=\FSF(G).
	\end{equation}
	In particular, the FSF on each component of $\ol{\WSFs(G)}$ and $\ol{\FSFs(G)}$ has a unique component a.s.  
\end{theorem}
An immediate corollary of Theorem~\ref{thm:resample} is 
\begin{corollary}\label{cor:wired1}
	$\WSF(\ol{\WSFs(G)})=\WSF(G)$  if and only if each component of $\ol{\WSFs(G)}$ has the property that $\FSF=\WSF$. \qed
\end{corollary}

Corollary~\ref{cor:wired1} implies Theorem~\ref{thm:resample-Z} as follows.  Recall $\cE$ defined as  in \eqref{eq:energy}.
A necessary and sufficient condition  for $\FSF(G)=\WSF(G)$ is  that
the only harmonic functions $f$ on $G$ with $\cE(f)<\infty$ are constant functions \cite{USF}. This is known to be the case when $G$ is transitive and \notion{amenable}, i.e.,
$\inf_{K} \#\partial K/\#K= 0$, where the infimum is over all finite vertex sets $K$ of $G$.
By \cite[Theorem 5.5]{BLS99}, every amenable transitive graph has the property that each component of every random subgraph with automorphism-invariant law also a.s.\ has no nonconstant harmonic functions $f$ with $\cE(f)<\infty$.
This gives Theorem~\ref{thm:resample-Z}. 

\begin{corollary}\label{cor:wired2}
	$\WSF(\ol{\WSFs(G)})=\WSF(G)$ (resp., $\FSF(\ol{\WSFs(G)})=\FSF(G)$) if  each component of $\WSFs(G)$ (resp., $\FSFs(G)$) is one-ended, that is, has a single end a.s.
\end{corollary} 
In \cite{End}, it is proved that  the one-end property of WSF
components holds for all transient vertex-transitive graphs (also see \cite[Theorem 10.49]{TreeBook}).  Thus
$\WSF(\ol{\WSFs(G)})=\WSF(G)$ in this case. This  in particular gives another proof of Theorem~\ref{thm:resample-Z}.
For more general results on the one-ended property of $\FSF$ and $\WSF$, see
\cite{End,Tom-End}.

Inspired by Morris' aforementioned  result that each component of $\WSF$ on every graph is a.s.\ recurrent, we conjecture that  	$\WSF(\ol{\WSFs(G)})=\WSF(G)$  for every locally finite, connected graph $G$ as in Theorem~\ref{thm:resample}.

Neither  $\WSF(\ol{\FSFs(G)}) =\FSF(G)$ nor $\WSF(\ol{\FSFs(G)}) =\WSF(G)$ holds for all graphs.
For counterexamples to the first equality, let $G$ be a tree with the property that $\WSFs(G)\neq G$ a.s. (For example, $G$ could be a regular tree.)  Then $\FSFs(G)=G$ a.s.\ while $\WSF(\ol{\FSFs(G)})= \WSF(G)$.  A counterexample for the second equality will be given in Section~\ref{subsec:counter}. 

\bigskip
\noindent{\bf Acknowledgment.}
We thank Pengfei Tang for assistance in completing the proof of Proposition \ref {prop:loop}, which simplifies arguments in an earlier version of our paper.
We are grateful to the referee for several comments that improved the paper.
This work was begun while the third author was an intern in the Theory Group at Microsoft Research, Redmond.
\section{Preliminaries}\label{sec:pre}
\subsection{Basic notations}\label{subsec:notion}
The set of positive integers is denoted by $\N$. Given a finite set $A$, we write $\# A$ for the cardinality of $A$. 
Given two sets $A,B$, their symmetric difference $(A\setminus B)\cup (B\setminus A)$ is denoted by $A\xor B$.
We use the asymptotic notation that two nonnegative
functions $f(x)$ and $g(x)$ satisfy $f\lesssim g$ if
there exists a constant $C>0$ independent of $x$  such that $f(x)\le Cg(x)$. We write $f\gtrsim g$ if $g\lesssim f$ and write $f\asymp g$ if $f\lesssim g$ and $f\gtrsim g$. 

Given a graph $G$, write $V(G)$ and $E(G)$ for the vertex and edge sets of $G$, respectively. 
When $G=\Z^d$ for some $d\in\N$, we write $o$ for its origin. If $v,u\in V(G)$ are adjacent, we write $v\sim u$ and write $(u,v)$ for the edge between them.
For $v\in V(G)$, let $\deg(v)$ be the degree of $v$, 
which is the number of vertices adjacent to $v$.
In our paper, graphs are assumed to  be \notion{locally finite}, that is, $\deg(v)<\infty$ for every $v\in V(G)$.
A graph $H$ is called a \notion{subgraph} of $G$ if $V(H)\subset V(G)$ and $E(H)\subset E(G)$. 
If $H$ and $H'$ are subgraphs of $G$, we write $E(G)\setminus E(H)$ as $G\setminus H$ and
$E(H)\xor E(H')$ as $H\xor H'$.

Given a family of probability measures $\Seq{\mu_t}_{t\in T}$ with index set
$T$, a \notion{coupling} of $\Seq{\mu_t}_{t\in T}$ is a family of random variables
$\Seq{X_t}_{t\in T}$ on a common probability  space such that $X_t$ is
distributed as $\mu_t$ for all $t\in T$.  Suppose $A$ and $B$ are two
probability measures on  the space of subgraphs of a graph $G$.  If there is a coupling $(\mathfrak
A,\mathfrak B)$ of $(A,B)$ such that $\mathfrak A\subset \mathfrak B$ a.s., we say that
$A$ is \notion{stochastically dominated} by $B$, written as $A \preccurlyeq B$. 

Let $I$ be an interval in $\Z$. 
Suppose $\path=\Seq{v_i}_{i\in I}$  is a sequence of vertices in $G$ indexed by $I$ such that $v_i\sim v_{i+1}$ whenever $i$ and $i+1$ are both in $I$. 
Then we call $\cP$ a   \notion{path} in $G$.
If $v_i\neq v_j$ as long as $i\ne j$, we say $\cP$ is \notion{simple}. 
If $I=\{0,\ldots, n\}$, then $\cP$ is called a finite path and  $|\cP|\defeq n$ is called the  \notion{length} of $\cP$.
We call the path $\Seq{v_{n-i}}_{0\le i\le n}$ the \notion{reversal} of $\cP$.
If we further have $v_0=v_n$, then  we call $\path$   a (rooted) \notion{loop}\footnote{This is a topological loop, also called a cycle in graph theory,	as opposed to the term ``loop" in graph theory.}
and $v_0$ the \notion{root} of $\cP$. 
If $I=\N\cup\{0\}$ (resp., $I=\Z$), we call $\cP$ an \notion{infinite} (resp., \notion{bi-infinite}) path. 
We call $t$  a \notion{cut time} of $\cP$ if $\{v_i\}_{i<t} \cap \{v_i\}_{i>t}=\varnothing$.

Given $x,y\in V(G)$, let $d_G(x,y)$ be the minimal length of a path starting from $x$ and ending at $y$ if $x,y$ are in the same component of $G$ and $\infty$ otherwise. 
We call $d_G(\cdot,\cdot)$ the \notion{graph distance} on $G$.
For $v\in V(G)$ and $r>0$, let $B_G(v,r):=\{x\in V(G) \st d_G(v,x)\le r\}$. We identify $B_G(v,r)$ with its induced subgraph.

A  graph is called a \notion{forest} if for any pair of distinct vertices there exists at most one simple path connecting them. A connected forest is called a \notion{tree}. 
Given a  connected graph  $G$, a \notion{spanning tree} (resp., \notion{forest}) on $G$ is a subgraph $T\subset G$ such that $T$ is a tree (resp., forest) and $V(T)=V(G)$.

The simple random walk on $G$ is the Markov chain $\Seq{S(n)}_{n\ge 0}$  on the state space $V(G)$ such that 
$\P[S(n+1)=u \mid S(n)=v]=\deg (v)^{-1}$ for all $u\sim v$ and $n\ge 0$.
The \notion{transition kernel} $p$ of $G$ is defined by $p_t(x,y)=\P[S(t)=y]$ for $x,y\in V(G)$ and $ t\in \N\cup\{0\}$, where $S$ is a simple random walk on $G$ starting from $x$.
When $G=\Z^d$, it is well known  that $p_t(o,o)\asymp t^{-d/2}$ for even $t \ge 2$.

\subsection{Wilson's algorithm}\label{subsec:Wilson}
Given a finite path $\path=\Seq{v_i}_{0\le i\le n}$ in a graph $G$ of length $n\in\N$,  
the \notion{(forward) loop erasure} of $\path$ (denoted by $\LE[\path]$) is the path defined
by erasing cycles in $\path $ chronologically. More precisely,  we define $\LE[\path]$ inductively as follows. The first vertex $u_0$ of $\LE[\path]$ equals $v_0$.
Supposing that $u_j$ has been set, let $k$ be the last index such that $v_k=u_j$. Set $u_{j+1}\defeq v_{k+1}$ if $k<n$; otherwise, let
$\LE[\path]=\Seq{u_i}_{0\le i\le j}$.  If $\path$ is an infinite path that visits no vertex infinitely many times, then we define $\LE[\path]$ in a similar fashion. In particular,
if $S$ is a sample of  simple random
walk on a transient graph $G$, then $\LE[S]$ is defined a.s. In such a case, we
call the law of $\LE[S]$ the \notion{loop-erased random walk} (LERW) on $G$. 

In \cite{Wilson}, Wilson discovered an algorithm for sampling uniform spanning trees on finite graphs using loop-erased random walk. In
\cite{USF}, Wilson's algorithm was adapted to sample $\WSF$ on a
transient graph $G$. This method is called 
\notion{Wilson's algorithm rooted at
infinity}, which we now review. The algorithm goes by sampling a
growing sequence of subgraphs of $G$ as follows. Set $\cT_0\defeq \varnothing$.
Inductively, for each $n\in \N$, choose $v_n\in V(G)\setminus V(\cT_{n-1})$ and run a simple random walk starting at $v_n$. 
Stop the walk when it hits $\mathcal{T}_{n-1}$ if it does; otherwise, let it run
indefinitely. Denote the resulting path by $\path_n$, and set $\mathcal T_n\defeq \mathcal T_{n-1}\cup \LE[\path_n]$. Write $\WSFs := \bigcup_n \mathcal{T}_n$. According to \cite[Theorem
5.1]{USF},  no matter how $\Seq{v_n}_{n\ge 1}$ are chosen, as long as $V(\WSFs) = V(G)$, the law of 
$\WSFs$ is $\WSF(G)$.

\subsection{Bounds on effective resistance}\label{sub:eff}
 Nash-Williams' inequality (see, e.g., \cite[Section 2.5]{TreeBook}) is a useful lower bound for the effective resistance.
Here we record a generalization of  Nash-Williams' inequality.
\begin{lemma}\label{lem:NW}
Given a graph $G$ with two disjoint subsets $A$ and $B$ of $V(G)$, a set $\cC\subset E(G)$ is called a  \notion{cut set} between $A$ and $B$ if\/ $\forall o\in A$ and $\forall z\in B$, every path from $o$ to $z$  must use an edge in $\cC$.
Suppose $\cC_1, \ldots, \cC_n$ 
are cut sets between $A$ and $B$ for some $n\in \N$.
For $e\in E$, let $j(e) \defeq  \#\{k  \st e \in \cC_k\}$.
Then   $\eff^{G}(A,B)\ge\sum_{k=1}^n \Bigl(\sum_{e \in \cC_k}
j(e) c(e)\Bigr)^{-1}$.
\end{lemma}
\begin{proof}
The proof is the same as the classical case in \cite[Section 2.5]{TreeBook}, 
with a slight modification when applying the Cauchy--Schwarz inequality.
We leave the details to the reader.
\end{proof}
This is equivalent to having disjoint cutsets in a subdivided network. 
In that form but for recurrence of infinite, locally finite networks, it was given by \cite{NashWil:crit},
who also proved that the existence of such cutsets in a subdivided network
is necessary for recurrence.
See \cite{McGuinness:NW} for an extension to non-locally finite networks.

The next lemma says that effective  resistance is stable under local modification.
\begin{lemma}\label{lem:eff}
Suppose $H$ and $H'$ are two connected subgraphs of a graph $G$ such that $\#(H\xor H')<\infty$.
Then there exists a constant $c>0$ depending on $G,H,H'$ such that 
$\eff^{H}(u,v)\le \eff^{H'}(u,v)+c$ for all
$u,v\in V(H)\cap  V(H')$. 
\end{lemma}
\begin{proof}
It suffices to show that if $H'\subset H$ and $\#\bigl(E(H)\setminus E(H')\bigr)=1$, then there exists $c'>0$ depending only on $H$ and $H'$ but not  on $u,v\in V(H')$
such that 
\begin{equation}\label{eq:compare}
0\le \eff^{H'}(u,v)-\eff^{H}(u,v)\le  c'.
\end{equation} 
Once this is proved, a similar statement then follows for $H'\subset H$ and $\#\bigl(E(H)\setminus E(H')\bigr)< \infty$. Then the general case follows by comparing both $H$ and $H'$ to the union graph $H \cup H'$ of the two.

By Rayleigh's monotonicity principle (see, e.g., \cite[Section~2.4]{TreeBook}), adding an edge can only decrease the effective resistance, hence $\eff^{H'}(u,v)\ge\eff^{H}(u,v)$.
To prove the other direction of \eqref{eq:compare}, 
we use Thomson's principle (see, e.g., \cite[Section~2.4]{TreeBook}) that the effective resistance between two vertices is the minimum
energy (i.e., the sum of the squares of all edge flows) among all unit flows between the two vertices. 
We may start from the minimizing flow for $H$ from $u$ to $v$ and then construct a flow on $H'$ 
between the same vertices by replacing the current flow along the removed edge $e$ with a flow along a path in $H'$ connecting 
the two endpoints of $e$. This  increases the flow energy by an additive constant that depends only on $H$ and $H'$.
\end{proof}

\subsection{Indistinguishability  of WSF components}
In this subsection, we review a basic ergodic-theoretic property of components in $\WSF$ on transient vertex-transitive graphs.
We call a triple $(G,\rho,\omega)$ a \notion{subgraph-decorated rooted graph} if $G$ is a locally finite, connected graph, $\rho$ is a distinguished vertex in $G$ called the \notion{root}, and $\omega$ is a function from $E(G)$ to $\{0,1\}$. We think of $\omega$ as a distinguished subgraph spanned by the edges $\{e\in E(G) \st \omega(e)=1\}$. 
Given two such triples $(G,\rho,\omega)$  and $(G',\rho',\omega')$, an isomorphism between them is a graph isomorphism between $G$ and $G'$ that preserves the root and the subgraph.
Write $\omega_{\rho, r}$ for the restriction of $\omega$ to $B_G(\rho, r)$.
Let $\cG^{\{0,1\}}_\bullet$ be the space  of subgraph-decorated rooted graphs modulo isomorphisms. 
We  endow $\cG^{\{0,1\}}_\bullet$ with the \notion{local topology} where two elements $(G,\rho,\omega)$ and $(G',\rho',\omega')$ in $\cG^{\{0,1\}}_\bullet$ are close 
if and only if $(B_G(\rho,r),\rho,\omega_{\rho, r})$ and $(B_{G'}(\rho',r),\rho',\omega'_{\rho', r})$ are isomorphic to each other for some large $r$. 

Given $(G,v,\omega)\in\cG^{\{0,1\}}_\bullet$, we define $K_\omega(v)$ to be the connected component of $v$ in $\omega$.
A Borel-measurable set $\cA \subset \cG^{\{0,1\}}_\bullet$ is called a \notion{component property} if 
$(G, v, \omega) \in \cA $ implies $(G, u, \omega) \in \cA$ for all $u \in K_\omega(v)$.
Given a component property $\cA$, we say that a connected component $K$ of $\omega$ has property $\cA$ if $(G,u,\omega)\in \cA$ for some (and  equivalently every) $u\in V(K)$.
A component property $\cA$ is called a \notion{tail component property} if $(G,v,\omega)\in \cA$ implies $(G,v,\omega')\in \cA$ 
for all  $\omega'\subseteq E(G) \text{ such that } \omega \xor \omega' $ and $K_\omega(v)\xor K_{\omega'}(v)$  are both finite.

As a corollary of  \cite[Theorem 1.20]{Indistinguish}, we have
\begin{lemma}
\label{lem:0-1}
Suppose $G$ is a transient, vertex-transitive graph. For every tail component property $\cA$,
either almost surely every connected component of\/ $\WSF(G)$ has property $\cA$, or almost surely  none of the connected components of\/ $\WSF(G)$ have property $\cA$.
\end{lemma} 
By Lemma~\ref{lem:eff}, for a vertex-transitive graph, both the properties in Theorems~\ref{thm:recurrence} and~\ref{thm:linear} are tail component properties.
Therefore, we have
\begin{lemma}\label{lem:compoent}
Consider $\WSF(\Z^d)$ for $d\ge 5$. Recall the notations $\cT_v$ and $\ray_v$ in Definition~\ref{def:rays}. Let $E_v$ be the event that $\ol\cT_v$ is recurrent
and $F_v$ be the event $\liminf_{n\to \infty}n^{-1}\eff^{\ol{\cT_v}}\bigl(v, \ray_v(n)\bigr) >0$. Then neither $\P[E_v]$ nor $\P[F_v]$ depends on $v$. Moreover, 
both $\P[E_v]$ and $\P[F_v]$ belong to $\{0,1\}$.
The same holds with $\Z^d$ replaced by any transient, vertex-transitive graph.		\qed
\end{lemma}

\subsection{Two-sided random walk and loop-erased random walk}\label{subsec:two-sided}
For $d\in \N$, let $S^1$ and $S^2$ be two independent simple random walks on $\Z^d$ starting from the origin of $\Z^d$.
For $n\in \Z$, let $S(n):=S_1(n)$ if $n\ge 0$ and $S(n):=S^2(-n)$ if $n<0$. 
We call the law of the bi-infinite path $\Seq{S(n)}_{n\in\Z}$ the \notion{two-sided  random walk} on $\Z^d$.
It is standard that $\P[S_1([0,\infty)) \cap S_2([1,\infty))=\varnothing]>0$ if and only if $d\ge 5$ (see, e.g., \cite[Theorem~10.24]{TreeBook}).
For $d\ge 5$, consider the event 
\begin{equation}\label{eq:B}
E := \bigl\{\LE\big[S^1\big](m)\neq S^2(n) \textrm{ for all } m\ge 0,\,n\ge 1\bigr\}.
\end{equation}
Since with positive probability $0$ is a cut time of $S$, we have that   $\P[E]>0$.
Define $\wt S(n)$ to be $\LE\big[S^{2}\big](-n)$
for $n\le 0$ and $\LE\big[S^{1}\big](n)$ for $n\ge 0$.
The conditional law of 	$\Seq{\wt S(n)}_{n\in \Z}$ conditioned on $E$ is called the \notion{two-sided loop-erased random walk} on $\Z^d$. 
It is clear that without loop-erasures, $\Seq{S(n+1)-S(n)}_{n\in \Z}$ is stationary and ergodic; indeed, it is an IID sequence. In fact, two-sided LERW also has stationary, ergodic increments:
\begin{lemma}\label{lem:bi}
Suppose $X$ is two-sided loop-erased random walk on $\Z^d$ for $d\ge 5$.
Then $\Seq{X(n+1)-X(n)}_{n\in \Z}$ is stationary and ergodic.
\end{lemma}
Lawler \cite{Lawler-Bi} introduced the two-sided LERW on $\Z^d$ ($d\ge 5$) and showed that it is the local limit of the usual LERW  viewed from nodes with large index.
An essential ingredient to the proof of Lemma~\ref{lem:bi} is the reversibility of the loop-erasing operation for simple random walk, which was also first proved in \cite{Lawler-Bi}. 
Given the reversibility, we observe that Lemma~\ref{lem:bi} can be deduced from the ergodicity of the two-sided random walk and 
the following basic fact from ergodic theory (see, e.g., \cite{ergodic}).
\begin{lemma}[Kac's Lemma]\label{lem:Kac}
Suppose $\Omega$ is a measurable space and $T\colon\Omega\to\Omega$ is measurable.
Suppose $\P$ is a probability measure on $\Omega$ which is preserved by $T$ and is ergodic. Let $A\subset \Omega$ be an event such that $\P[A]>0$ and 
let $\tau(\omega):=\inf \{n\in \N \st T^n(\omega)\in A\}$ for all $\omega\in\Omega$. 
Let $T_A(\omega):=T^{\tau(\omega)}(\omega)$ for all $\omega\in A$.
Then $T_A$ is an ergodic, measure-preserving map from $A$ to $A$ under the conditional probability measure $\P[\cdot\mid A]$. Moreover, $\E[\tau\mid A] =\P[A]^{-1}$.
\end{lemma}
To put Lemma~\ref{lem:bi} in the setting of Lemma~\ref{lem:Kac}, let us consider the two-sided simple random walk $S$. Since $S$ can be almost surely decomposed into finite paths separated by cut times, the \emph{forward} loop-erasure of the  path $\Seq{S(n)}_{n\le 0}$  is well defined, which we denote by $\LE[S(-\infty,0]]$. 
By the reversibility of the loop-erasing operation,  the path $\LE[S(-\infty,0]]$ has the same law as $\LE[S[0,\infty)]$.
Now we use the event $A\defeq \{ \LE[S(-\infty,0]] \cap S[1,\infty) =\emptyset \}$, which plays the same role as the event $E$ in $\eqref{eq:B}$. Let $T$ be the forward shift operator of $\Seq{S(n)}_{n\in \Z}$. Applying Lemma~\ref{lem:Kac} to $T$ and $A$, we get Lemma~\ref{lem:bi}. 

Using estimates for random walk on $\Z^d$, it was shown in \cite{Lawler-Bi} that the two-sided LERW is weakly mixing, which is a property stronger than ergodicity. For this paper, we need  only stationarity (see Section~\ref{sec:highD} for its use) and our argument  can be readily extended to more general unimodular vertex-transitive graphs.

A network is a graph with positive weights on its edges. The corresponding network random walk is the Markov chain on its vertices that moves to a neighboring vertex of its present state with probabilities proportional to the weights of the incident edges. Constant weights correspond to simple random walk.
For a vertex $x$ in a general transient network, $G$, define 
\begin{align*}
Z_i(x):=\sum_{t = 0}^{\infty}(t+1)^i p_t(x, x).  
\end{align*}  
When $G$ is $\Z^d$, or, more generally, transitive, $Z_i(x)$ does not depend on $x$, and so we will write simply $Z_i$ in the transitive case.
Because $Z_0(x)$ is the expected number of visits to $x$ by a network random walk $S$ starting from $x$, we have that $Z_0(x) = 1/\P\bigl[\forall n \ge 1\enspace S(n) \ne x\bigr]$. If $L(x)$ denotes the time of the last visit to $x$, it follows that for any loop $\gamma$ rooted at $x$, we have 
\begin{equation}\label{eq:normalize}
\P\bigl[S\bigl([0, |\gamma|]\bigr) = \gamma\bigr]/Z_0(x) = \P\bigl[S\bigl([0, L(x)]\bigr) = \gamma\bigr],
\end{equation}
and so
\begin{equation}\label{eq:normalizeExp}
\E\bigl[L(x) + 1\bigr] = Z_1(x)/Z_0(x).
\end{equation}

\begin{lemma}\label{lem:cond}
Let $G$ be a network and $x$ be a vertex of $G$. Let $A$ be a set of
vertices such that there is positive probability that the network walk $S$ from
$x$ avoids $A$ forever: $\P[\forall n \ge 0 \enspace S(n) \notin A] > 0$.  Let $L \defeq \sup \{n \ge 0 \st S(n) = x\}$ be the last time the walk is at $x$.  
Then the probability that
the first $L$ steps of the walk is equal to a particular loop $\gamma$ given that
$A$ is avoided forever is at most the unconditional probability that the
first $|\gamma|$ steps of the random walk is $\gamma$:
$$
\P\bigl[S\bigl([0, L]\bigr) = \gamma \bigm| \forall n \ge 0 \enspace S(n) \notin A\bigr]
\le
\P\bigl[S\bigl([0, |\gamma|]\bigr) = \gamma].
$$
\end{lemma}

\begin{proof}
This is a simple modification of the proof of the display after (5.4) of \cite{Hutchcroft:universality}.
To be more precise, given  a path $\gamma$ avoiding $A$, we have
\begin{align*}
\P\bigl[S\bigl([0, L]\bigr) = \gamma \bigm| \forall n \ge 0 \enspace S(n) \notin A\bigr]
&=
\frac{\P\bigl[S\bigl([0, |\gamma|]\bigr) = \gamma,\  \forall n > |\gamma| \enspace S(n) \notin A \cup \{x\}\bigr]}
     {\P\bigl[\forall n \ge 0 \enspace S(n) \notin A\bigr]}
\\ &=
\frac{\P\bigl[S\bigl([0, |\gamma|]\bigr) = \gamma] \cdot \P[\forall n > |\gamma| \enspace S(n) \notin A \cup \{x\}\bigr]}
     {\P\bigl[\forall n \ge 0 \enspace S(n) \notin A\bigr]}
\\ &\le
\P\bigl[S\bigl([0, |\gamma|]\bigr) = \gamma].
\qedhere
\end{align*}
\end{proof}

\begin{proposition}\label{prop:loop}
Let $G$ be a transitive network and $o$ be a vertex of $G$. Consider the network walk $S$ from $o$.
Define \(K_i\defeq  \#\{s \st \big|\LE[S([0,s])]\big |=i \}\).
Then $\E [K_i]\le Z_1$ for all $i \ge 0$.
\end{proposition}
\begin{proof}
That $\E[K_0] = Z_0 \le Z_1$ is trivial.
Define $T_0 := 0$ and for $i\in \N$, let 
\[
T_i\defeq 1 + \sup \bigl\{t \ge T_{i-1} \st S(t) = S(T_{n-1})  \bigr\}.
\]
Then $\big|\LE([S[0, T_j])\big| = j$ for all $j \ge 0$.
Defining \(K_{i, j} \defeq  \#\{s \in [T_j, T_{j+1}) \st \big|\LE[S([0,s])]\big| =i \}\) for $0 \le j \le i$, we have 
$K_i = \sum_{j=0}^{i} K_{i, j}$. 

Let  $S'$ be an independent simple random walk from $o$;
let $T'$ be one more than the last time $S'$ visits $o$. Let $R_j :=
\{s \in [0, T') \st \big|\LE[S'([0,s])]\big| = i - j \}$; these sets are
pairwise disjoint.  By Lemma~\ref {lem:cond} and equation \eqref {eq:normalize}, for $i \ge j \ge 0$ and $k \ge 0$, 
\[
\P\bigl[K_{i, j} = k \bigm|  S[0, T_j]\bigr] 
\le
Z_0 \,\P[\# R_j = k],
\]
whence 
\[
\E[K_i]
=
\sum_{j=0}^{i} \E[K_{i, j}]
\le
\sum_{j=0}^{i} Z_0 \,\E[\# R_j]
=
Z_0 \,\E\Bigl[\#\bigcup_{j=0}^i R_j\Bigr] \le Z_0 \,\E[T']
=
Z_1,
\]
where, in the last step, we used \eqref {eq:normalizeExp}. 
\end{proof}

\section{Recurrence when \texorpdfstring{$d\ge 8$}{d >= 8}}\label{sec:8D}
In this section, we first prove Theorem~\ref{thm:recurrence} and then extend the result to vertex-transitive graphs in Section~\ref{subsec:transitive}.
Recall $\cT_o$ and $\ray_o$ as in Definition~\ref{def:rays} for the origin $o$ of $\Z^d$.
Given $n\in \N\cup\{0\}$, we call the connected component of $\tree \setminus \bigl(\ray_o[0,n-1]\cup\ray_o[n+1,\infty)\bigr)$ containing $\ray_o(n)$ the $n$th \notion{bush} of $\cT_o$ and denote it by $\bush_n$.
Given an edge $e$ in $\Z^d$ and two subgraphs $H_1$ and $H_2$ of $\Z^d$ with $V(H_1)\cap V(H_2)=\varnothing$,
we say that $e$ \notion{joins} $H_1$ and $H_2$ if one endpoint of $e$ is in  $H_1$ and the other is in $H_2$.
\begin{lemma}\label{lem:expectation}
For $0\le j\le n$ and $\ell\ge 1$, let $N_{j,\ell}(n)$ be the number of edges  joining  $\bush_{n-j}$ and $\bush_{n+\ell}$.	
Suppose $d\ge 8$. Then there exists a constant $C>0$ such that 
\begin{equation}\label{eq:Onm}
\sum_{0\le j\le n} \sum_{\ell \ge m} \E[N_{j,\ell}(n)] \le  C \log \Bigl(1 + \frac {n}{m}\Bigr) \qquad \textrm{for all } m,n\in \N.
\end{equation}
\end{lemma}
We postpone the proof of Lemma~\ref{lem:expectation} to Section~\ref{subsec:lem} and proceed to prove Theorem~\ref{thm:recurrence}.
\begin{proof}[Proof of Theorem~\ref{thm:recurrence}]
By Lemma~\ref{lem:compoent}, we see that Theorem~\ref{thm:recurrence} is equivalent to the statement  that $\graph$ is recurrent a.s.

By \eqref{eq:Onm} with $m=1$, there exists a constant $C>0$ and a sequence  $n_k\in [k^{2k},2k^{2k}]$ such that 
\begin{equation}\label{eq:log2}
\sum_{0\le j\le n_k} \sum_{\ell \ge 1} \E[N_{j,\ell}(n_k)]\le C\log n_k\qquad\textrm{for all }k\in\N.
\end{equation}
Define $\cC_k$ to be the set consisting of
\vadjust{\kern2pt}%
edges  joining $\bigcup_{m\le n_k} \bush_m $ and  $\bigcup_{m > n_k} \bush_m $.  Then removing $\cC_k$ from $\graph$ leaves $o$ in a finite component.
Since $\E[\#\cC_k]\lesssim \log n_k$ by \eqref{eq:log2} and
\[
\sum_{k=1}^{\tau}\log\left(n_{k}\right)\lesssim \tau^2\log \tau\qquad \textrm{for all } \tau \ge 1,
\]
the argument in \cite[Lemma 13.5 and Remark 13.6]{USF} yields \(\sum_{1}^{\infty}(\#\cC_k)^{-1}=\infty\) a.s.

Let $I_k$  be  the event that there exists an edge joining $\bigcup_{m\le n_k} \bush_m $ and  $\bigcup_{m \ge n_{k+1}} \bush_m $. 
Since \(\sum_{1}^{\infty}\frac{n_k}{n_{k+1}-n_k}<\infty,\)
by \eqref{eq:Onm} and the Borel--Cantelli   lemma, we know that almost surely  only finitely many events $I_k$ occur.
Therefore there exists a (random) \vadjust{\kern2pt}%
$K\in\N$ such that the elements in $\{\cC_k \st k\ge K\}$ are all disjoint.
By the Nash-Williams criterion (see, e.g. \cite[Sec.~2.5]{TreeBook}), it follows that $\graph$ is recurrent a.s.
\end{proof}
\subsection{Proof of Lemma~\ref{lem:expectation}}\label{subsec:lem}

Fix $d\in \N$. For $z\in \Z^d$, let $\Omega^z$ be the space of loops in $\Z^d$ rooted at $z$ (see Section~\ref{subsec:notion} for the definition). 
Define a measure  $\mu$ on $\Omega^z$ by requiring  $\mu(\gamma)\defeq (2d)^{-|\gamma|}$ for all $\gamma\in\Omega^z$.
We call $\mu$ the \notion{loop measure} and $\mu(\gamma)$ the \notion{weight} of $\gamma$.
Here we drop the dependence of $\mu$ on $z$ for simplicity of notation.
In different places, we will consider loops with additional markings. For example, 
let $\Omega^z_{\op s}\defeq \{(\gamma,i) \st \gamma\in \Omega^z,\; i\in [0,|\gamma|-1]\cap \Z \}$. 
Each element in $\Omega^z_{\op s}$ is a loop $\gamma$ rooted at $z$ with a \notion{marked step} being the ordered pair $\Seq{\gamma(i),\gamma(i+1)}$.
By assigning  each element in $\Omega^z_{\op s}$ the weight of its loop, we  define a measure on $\Omega^z_{\op s}$, which we still denote by $\mu$ in a slight abuse of notation.

\begin{proof}[Proof of Lemma~\ref{lem:expectation}]
By linearity of expectation, we estimate  $\E[N_{j,\ell}(n)]$ by estimating the probability of joining $\bush_{n-j}$ and $\bush_{n+\ell}$  for each edge of $\Z^d$.  
Let  $x$ and $y$ be two adjacent vertices in $\Z^d$ and $S$, $S^1$, and $S^2$ be three independent simple random walks on $\Z^d$ starting from $o$, $x$, and $y$ respectively. 
Suppose that $\WSF(\Z^d)$ is sampled via Wilson's algorithm rooted at infinity by first sampling $S$, $S^1$, and $S^2$, and then other random walks. 
Fix $0\le j\le  n $ and $\ell, m \in  \N$. 
Given $s,s',t'\in\N\cup\{0\}$  and $z,w\in \Z^d$,   let $E^{x,y}(s,s',t',z,w)$ be the event that
	\begin{enumerate}
		\item $S(s)=z $ and $\big |\LE[S([0,s])]\big |=n-j$;
		\item $S^1(s')=z$ and $S^2(t')=w$;
		\item $\lambda\defeq \sup\{k \st S(k)=w\}\in [s, \infty)$ and $\big |\LE\big [S([s,\lambda])\big ]\big |=j+\ell$; and
		\item $t'=\inf\left\{k \st S^2(k)\in  \LE\big[S([s,\lambda])\big]\right\}$.
	\end{enumerate}
Then $\bigl\{ x\in V(\bush_{n-j})\textrm{ and }y\in V(\bush_{n+\ell})\bigr\}\subset \bigcup E^{x,y}(s,s',t',z,w)$, where the union ranges over all possible tuples $(s,s',t',z,w)$.

Recall $\Omega^z_{\op s}$ defined right above.  Let $\Omega^z_{\op s}(x,y,w,s',t')\subset \Omega^z_{\op s}$ be the set of $(\gamma,i)$ that satisfy
\begin{enumerate}
\item $\gamma(i)=y$ and $\gamma(i+1)=x$;
\item  $s'=|\gamma|-i-1$;
\item $\LE[\gamma([0,i])](j+\ell)=w$; and
\item $\max\{k\le i \st \gamma(k)=w\}=i-t'$.
\end{enumerate}
On $E^{x,y}(s,s',t',z,w)$, by concatenating  $S([s,\lambda])$, the reversal of $S^2([0,t'])$, the edge from $y$ to $x$, and $S^1([0,s'])$, 
we obtain an element in $\Omega^z_{\op s}(x, y, w, s', t')$ whose marked step is $\Seq{y,x}$. Therefore $\P[E^{x,y}(s,s',t',z,w)]$ equals 
\begin{equation}\label{eq:path}
\P\bigl[S(s)=z \textrm{ and }\big |\LE[S([0,s])]\big |=n-j\bigr]\cdot (2d) \cdot \mu[\Omega^z_{\op s}(x,y,w,s',t')],
\end{equation}
where the factor $2d$ comes from the fact that the step $\Seq{y,x}$ need not be traversed by $S$, $S^1$, or $S^2$.
Note that \(\Omega^z_{\op s}(x,y, w,s',t') \subset  \{(\gamma,i)\in \Omega^z_{\op s} \st |\gamma|\ge j+\ell\}\).
Now let $x$ and $y$ vary.
For different tuples $(x,y, w,s',t')$, the corresponding sets $\Omega^z_{\op s}(x,y, w,s',t')$  are disjoint (because the definition (1)--(4) determines $x$, $y$, $w$, $s'$ and $t'$ from $(\gamma, i)$).
Therefore 
\begin{equation}\label{eq:inclusion}
\sum_{x,y,w,s',t'}\mu[\Omega^z_{\op s}(x,y, w,s',t')] \le \mu[\{(\gamma,i)\in \Omega^z_{\op s} \st |\gamma|\ge j+\ell\}].
\end{equation}
Let $p$ be the transition kernel of $\Z^d$.
By the definition of  $\Omega^z_{\op s}$ and $\mu$, for all $t\in \N$, we have 
	\[ \mu\bigl[\{(\gamma,i)\in \Omega^z_{\op s} \st |\gamma|=t\}\bigr]=tp_t(z,z)=tp_t(o,o).\]
Since $d\ge 8$ and $p_t(o,o)\lesssim t^{-d/2}$,  we see from \eqref{eq:inclusion} that
\[
\sum_{x,y,w,s',t'}\mu[\Omega^z_{\op s}(x,y, w,s',t')] \le\sum_{t\ge j+\ell}	 tp_t(o,o) \lesssim (j+\ell)^{-2}.	
\]
Set 	\(K_i\defeq  \#\{s \st \big|\LE[S([0,s])]\big |=i \}\) for all $0\le i\le n$, as in Proposition~\ref{prop:loop}.	By \eqref{eq:path}, we have
\begin{multline*}
\sum_{x,y,s,s't',z,w}\P[E^{x,y}(s,s',t',z,w)] \lesssim (j+\ell)^{-2}\sum_{s,z}\P\bigl[S(s)=z, \big |\LE[S([0,s])]\big |=n-j\bigr]\\
=  (j+\ell)^{-2}\sum_{s}\P\bigl[\big |\LE[S([0,s])]\big |=n-j\bigr]=\E[K_{n-j}](j+\ell)^{-2}
{} \lesssim (j+\ell)^{-2}
\end{multline*}
by Proposition~\ref{prop:loop}.	Since $\E[N_{j,\ell}(n) ] \le 	\sum_{x,y,s,s't',z,w}\P[E^{x,y}(s,s',t',z,w)]$, we see that
\begin{equation*}
	\sum_{0\le j\le n} \sum_{\ell \ge m} \E[N_{j,\ell}(n)]  {} \lesssim  \sum_{0\le j\le n} \sum_{\ell \ge m} (j+\ell)^{-2}\lesssim \sum_{0\le j\le n} (j+m)^{-1} \lesssim \log \bigl(1 + n/m\bigr).
\qedhere
\end{equation*}
\end{proof}

\subsection{Extension to vertex-transitive graphs}\label{subsec:transitive}
Suppose $G$ is a vertex-transitive graph and $o\in V(G)$. We call $V(r)\defeq \#B_G(o,r)$ the \notion{volume growth function} of $G$. 
In this subsection, we explain the following extension of Theorem~\ref{thm:recurrence}.
\begin{theorem}\label{thm:recurrence1}
If $G$ is a  vertex-transitive graph with $V(r)\gtrsim r^8$, then almost surely each connected component of\/ $\ol{\WSF(G)}$ is recurrent.
\end{theorem}
Let $p$ be the transition kernel of $G$. 
It is standard that $p_t(o,o)\lesssim t^{-4}$ when $V(r)\gtrsim r^8$ (see, e.g., \cite[Corollary~6.32]{TreeBook}).
Given this, Proposition~\ref{prop:loop} implies Lemma~\ref{lem:expectation},  hence Theorem~\ref{thm:recurrence} in the same way with $\Z^d$ replaced by $G$ in Theorem~\ref{thm:recurrence}.

If $G$ is nonunimodular, then Theorem~\ref{thm:recurrence1} is essentially already known (see the definition of unimodular preceding Theorem~\ref{thm:linear1}).  In fact, it is well known that  $G$ is
\notion{nonamenable} in this case, that is, $\inf_{K}\#\partial K/\#K >0$, where the infimum is over all finite vertex sets $K$ of $G$; see, e.g., \cite[Proposition 8.14]{TreeBook}. Therefore,
by \cite[Theorem 13.1]{USF},  we have $\E \bigl[\#\bigl(\graph\cap B_{G}(0,n)\bigr)\bigr]\asymp n^2$, where $\tree$ is the component
of $\WSF(G)$ containing $o$. 
Now   \cite[Lemma 13.5]{USF} yields that $\graph$ is a.s.\ recurrent. 
Since $G$ is transient,  Lemma~\ref{lem:compoent} concludes Theorem~\ref{thm:recurrence1} in the nonunimodular case.

\section{Linear growth of resistance when \texorpdfstring{$d\ge 9$}{d >= 9}}\label{sec:highD}
Recall the two-sided LERW defined in Section~\ref{subsec:two-sided}. Now we  define the two-sided $\WSF$.
\begin{definition}\label{def:bi}
Given $d\ge 5$, sample a random spanning forest $\WSFs^2(\Z^d)$ on $\Z^d$ as follows.
	\begin{enumerate}
		\item Sample a two-sided loop-erased random walk $\wt S$.
		\item Conditioning on $\wt S$, sample a $\WSF$ (denoted by $\WSFs$) on the graph obtained from $\Z^d$ by identifying the trace of $\wt S$ as a single vertex.
		\item Set $\WSFs^2(\Z^d)$ to be the union of $\WSFs$ and the trace of $\wt S$, where $\WSFs$ is viewed as a random subgraph of $\Z^d$.
	\end{enumerate}
We call the law of $\WSFs^2(\Z^d)$ the \notion{two-sided wired spanning forest} on $\Z^d$ and denote it by $\WSF^2(\Z^d)$.
\end{definition}
It is clear that $\WSF^2(\Z^d)$ can be sampled from a modified version of Wilson's algorithm rooted at infinity: first sample a two-sided LERW and treat it as the first walk in  Wilson's algorithm; then proceed as in the original  Wilson's algorithm to form a spanning forest on $\Z^d$. The stationary two-sided LERW on $\Z^d$ was extended to $d=4$ in \cite{LSW-4D} and to $d=2,3$ in \cite{Lawler-two-sided} by a limiting procedure. Therefore $\WSF^2(\Z^d)$ can be defined for all $d\in \N$. However, we will not need the lower-dimensional cases.

By Lemma~\ref{lem:bi},  as a subgraph-decorated rooted graph, $(\Z^d,o,\WSFs^2(\Z^d))$ 
is stationary under shifting along the trace of $\wt S$.  We will use this stationarity and the ergodic theorem to prove Theorem~\ref{thm:linear}. 
The following lemma will be needed.
\begin{lemma}\label{lem:join}
In the setting of Theorem~\ref{thm:linear}, for $v\in\Z^d$ such that $v\neq o$, 
let $\cN_v$ be  the number of edges joining\footnote{Recall the notion from Section~\ref{sec:8D}.} $\tree$ and $\cT_v$ if $\tree\neq \cT_v$ and be 0 otherwise. Then $\E[\cN_v]<\infty$ for $d \ge 9$. 
\end{lemma}
\begin{proof}
We follow an argument similar to that in Lemma~\ref{lem:expectation}.  
Given two neighboring vertices $x$ and $y$, let $\cI_{x,y}$ be the event that $x\in \tree$ and $y\in \cT_v$. 
Suppose that $S^o$, $S^v$, $S^x$, and $S^y$ are independent simple random walks on $\Z^d$ starting from $o$, $v$, $x$, and $y$, respectively. By Wilson's algorithm, 
\[
\P [\cI_{x,y}]\le \P [S^x([0,\infty)) \cap S^o([0,\infty))\neq \varnothing\;\;\textrm{and}\;\; S^y([0,\infty)) \cap S^v([0,\infty))\neq \varnothing] . 
\]
Therefore
\begin{align}\label{eq:join-num}
\E [\cN_v] \le \sum_{x\sim y}\sum_{k,l,m,n \ge 0} \P
[S^o(k)=S^x(l)\;\;\textrm{and}\;\; S^y(m)=S^v(n)].
\end{align}
Now let $\Omega^{o,v}$ be the space of quadruples $(\gamma,\sigma,\tau,i)$ where
$\gamma$ is a path in $\Z^d$ from  $o$ to $v$, $\> \sigma,\tau\in [0,|\gamma|]\cap \Z$, and $i\in [0,|\gamma|-1]\cap \Z$. Here, $\sigma$ and $\tau$ are considered as  two marked times of $\gamma$ and $\Seq{\gamma(i),\gamma(i+1)}$ is considered as a marked step. 
Define the  measure $\mu$ on $\Omega^{o,v}$ by assigning weight $(2d)^{-|\gamma|}$ to each $(\gamma,\sigma,\tau,i)\in \Omega^{o,v}$. Then 
 \begin{equation*}
\mu[\Omega^{o,v}]=\sum_{t=0}^{\infty} \mu[\{(\gamma,\sigma,\tau,i)\in \Omega^{o,v} \st|\gamma|=t\}]=\sum_{t=0}^{\infty} t(t+1)^2p_t(o,v).
 \end{equation*}
 Since $p_t(o,v) \lesssim t^{-d/2}$ 
 (see \cite[Corollary~6.32(ii)]{TreeBook}), we see that $\mu[\Omega^{o,v}]<\infty$ if $d\ge 9$.

Let $\Omega^{o,v}(k,l,m,n)\defeq \{(\gamma, k,k+\ell+1+m, k+\ell) \in\Omega^{o,v} \st |\gamma|=k+\ell+1+m+n \}$.
By concatenating $S^o([0,k])$, the reversal of $S^x([0,\ell])$,  
the edge from $x$ to $y$, the path $S^y([0,m])$, and the reversal of $S^v([0,n])$, we see that
\(\sum_{x\sim y} \P [S^o(k)=S^x(l)\;\;\textrm{and}\;\; S^y(m)=S^v(n)]\) is no larger than \(\mu\bigl[\Omega^{o,v}(k,l,m,n)\bigr].\)
On the other hand, $\Omega^{o,v}(k,l,m,n)\cap\Omega^{o,v}(k',l',m',n')=\varnothing$ if $(k,l,m,n)\neq (k',l',m',n')$. Now interchanging the summations in \eqref{eq:join-num},  we get  $\E [\cN_v] \le\mu[\Omega^{o,v}]$, which is finite if $d\ge 9$.
\end{proof}
\begin{proof}[Proof of Theorem~\ref{thm:linear}]
By Lemma~\ref{lem:compoent}, it suffices to prove that
	\begin{equation}\label{eq:eff}
		\P\Bigl[\liminf_{n\to \infty}n^{-1}\eff^{\graph}(o,\ray_o(n))>0\Bigr] > 0.
	\end{equation} 	
	Let us make a particular choice in Wilson's algorithm rooted at infinity to sample $\WSF(G)$.
	\begin{enumerate}
		\item Sample a simple random walk $S^1$ from $o$ as the first walk in Wilson's algorithm.  
		\item Run an independent simple random walk $S^2$ from $o$. Define $\tau\defeq \inf\{t\ge 1 \st S^2(t)\notin  \LE \bigl[S^1\bigr]\}$ and $v\defeq S^2(\tau)$. 
		\item Use $S^2\bigl([\tau, \infty)\bigr)$ for the second simple random walk  in Wilson's algorithm.
		\item Sample the rest of $\WSF(G)$ according to Wilson's algorithm in an arbitrary way.
	\end{enumerate}	
Let $\P$ be the probability measure from the above sampling and let $\wt \P$ be $\P$ conditioned on the event $B\defeq\{\tau=1 \textrm{ and }\cT_v\neq \tree\}$.
Then
$B$ is exactly the event $E$ in~\eqref{eq:B}. 
We define $\wt S$ in terms of $(S^1,S^2)$ as in Lemma~\ref{lem:bi}, so that under $\wt\P$ it is a two-sided LERW.  
On the event $B$, let $\WSFs^2$ consist of the edges of $\WSF(G)$ and the edge $(o,v)$,  and let  $\wt\tree$ consist of the edges of $\cT_o$ and $\cT_v$ and the edge $(o,v)$.
By Lemma~\ref{lem:bi} and  Definition~\ref{def:bi},  
under $\wt \P$, we see that $\WSFs^2$ is distributed as  $\WSF^2(G)$ and $\wt\tree$  is the component of  $\WSFs^2$   containing $o$.

To prove \eqref{eq:eff}, recall the notion of $\bush_n$ in Section~\ref{sec:8D}. For $k\in\N$, let $\cC_k$
be the set of  edges  joining $\bigcup_{m\le k}\bush_m$ and  $\bigcup_{m\ge k+1}\bush_m$.
  For any edge $e$ of $\Z^d$, let $j(e)\defeq \#\{k \st e\in \cC_k \}$. 
Let $J_k:=\sum_{e\in \cC_k} j(e)$.  
Under $\wt \P$, for $n\in \Z$, let $\wt\bush_n$ be the connected
 component of $\wt\tree \setminus  \wt S(\Z\setminus \{n\})$ containing $\wt S(n)$.  
Let $\wt\cC_k$ be the set   of   edges  joining $\bigcup_{m\le k}\wt\bush_m$ and  $\bigcup_{m\ge k+1}\wt\bush_m$.
\vadjust{\kern2pt}%
Let  $\wt j(e)\defeq \#\{k \st e\in \wt \cC_k\}$ and  $\wt J_k:=\sum_{e\in \wt \cC_k} \wt j(e)$.
By Lemma~\ref{lem:join}, $\#\wt\cC_{-1}<\infty$  $\>\wt \P$-a.s.  
By the stationarity of  $\WSF^2(\Z^d)$, both $\Seq{\wt \cC_k}_{k\in\Z}$ and $\Seq{\wt J_k}_{k\in \Z}$ are stationary under $\wt \P$.
On the other hand, if $e\in \cC_k$ joins $\bush_m$ and $\bush_n$ for some $n>m$, we must have $ e\in \wt \cC_k$ and $j(e)=\wt j(e)=n-m$.
Therefore $\cC_k\subset \wt \cC_k$ and  $J_k\le \wt J_k<\infty$ $\>\wt \P$-a.s.\ for all $k\in\N\cup\{0\}$.
By the stationarity of  $\Seq{\wt J_k}_{k\in \Z}$ under $\wt \P$ and Birkhoff's ergodic theorem, there exists a random variable $Y$ such that  
\(\E^{\wt \P}[Y]=\E^{\wt \P}[\wt J_0^{-1}]>0\) and \(\lim_{n\to\infty} n^{-1}\sum^{n-1}_{k=0} \wt J_k^{-1}=Y\) $\>\wt \P$-a.s.
Since $J_k\le \wt J_k$, with positive probability under $\wt\P$ (hence under $\P$),  we have
\begin{equation*}
\liminf_{n\to\infty}n^{-1}\sum^{n-1}_{k=0} J_k^{-1}>0.
\end{equation*}
By Lemma~\ref{lem:NW}, $\eff^{\graph}(o,\ray_o(n))\ge \sum_{k=0}^{n-1} J_k^{-1}$, 
 which gives \eqref{eq:eff}.
\end{proof}
We conclude this section by the following straightforward extension of Theorem~\ref{thm:linear}.
Let $S$ be the two-sided random walk on $G$ defined  as in Section~\ref{subsec:two-sided} with $\Z^d$ and its origin replaced by $G$ and $o\in V(G)$.
A vertex-transitive graph $G$ is \notion{unimodular} if the automorphism group $\mathrm{Aut}(G)$ of $G$ is
unimodular, in other words,  $\mathrm{Aut}(G)$  admits a nontrivial Borel measure that is invariant under both left and right multiplication by group elements.
In this case, $G$ satisfies\footnote{It can be shown that unimodularity is, in fact, equivalent to \eqref{eq:unimodular}.}  
\begin{equation}\label{eq:unimodular}
\textrm{$\Seq{S(n)}_{n \in \Z}$ is stationary and ergodic viewed as path-decorated rooted graphs.}
\end{equation}
We will not  elaborate on the notion of unimodularity, but refer to \cite{BLPSgip} or \cite[Chapter 8]{TreeBook} for more background.
\begin{theorem}\label{thm:linear1}
Theorem~\ref{thm:linear} still holds if $\Z^d$ is replaced by a unimodular, vertex-transitive graph $G$ such that $V(r)\gtrsim r^9$.
(Recall the function $V(r)$ in Section~\ref{subsec:transitive}.) 
\end{theorem}
\begin{proof}
Note that $V(r)\gtrsim  r^9$ implies that the transition kernel satisfies $p_t(o,o)\lesssim t^{-9/2}$.
By inspection,  the proof of Theorem~\ref{thm:linear} still works given  this transition-kernel estimate and the fact that $\cT_o$ can be coupled with  the (stationary) two-sided $\WSF$ as in the proof of  Theorem~\ref{thm:linear} via Wilson's algorithm. 
This holds as long as  the two-sided LERW can be sampled from the  two-sided simple random walk as in Section~\ref{subsec:two-sided}. 
By \eqref{eq:unimodular} and Lemma~\ref{lem:Kac}, this is true if $G$ is unimodular.
\end{proof}
We expect that the unimodularity assumption in  Theorem~\ref{thm:linear1} can be removed. 
However, this would require a different approach, because for nonunimodular, vertex-transitive graphs, although the two-sided LERW can still  be defined  by a limiting procedure, it is not related to the two-sided simple random walk that we defined earlier.
\section{Resampling property}\label{sec:resample}
In this section, we  first prove Theorem~\ref{thm:resample} and Corollary~\ref{cor:wired2} in Section~\ref{subsec:resample}. 
Then we provide a counterexample to  $\WSF(\ol{\FSFs(G)})=\WSF(G)$ in Section~\ref{subsec:counter}.

\subsection{Proof of Theorem~\ref{thm:resample} and Corollary~\ref{cor:wired2}}\label{subsec:resample}
We introduce the following notation. Given a graph $G$, suppose $H$ is a random finite subgraph of $G$. Let us sample a random forest on $G$ as follows. First sample $H$. Conditioning on $H$, uniformly sample a spanning tree on each component of $H$. The unconditional law of the  resulting random forest is denoted by $\USF(H)$.
\begin{proof}[Proof of Theorem~\ref{thm:resample}]
We prove only $\FSF(\ol{\WSFs(G)})= \WSF(G)$ since $\FSF(\ol{\FSFs(G)})=\FSF(G)$ can be proved in exactly the same way. 

Fix $o\in V(G)$. For a positive integer $n$,  let $\WSFs^n$ be a sample of $\WSF(B_G(o,n))$. For $0 < m < n$,
thinking  of $\WSFs^n\cap B_G(o,m)$ and $\WSFs \cap B_G(o,m)$    as subgraphs of $B_G(o,m)$, let $K_{m,n}$ and $K_m$ be their induced-component graphs, respectively. 
For a fixed $m$, as $n$ tends to $\infty$, the laws of $\WSFs^n\cap B_G(o,m)$ and $\WSFs \cap B_G(o,m)$  can be coupled so that they are identical with probability $1-o_n(1)$.
Hence the same is true for $K_{m,n}$ and $K_m$. Conditioning on $K_{m,n} = K$, the conditional law of $\WSFs^n\cap B_G(o,m)$ is  $\USF(K)$ because every spanning forest of $K$ that is connected in each component of $K$ extends to a spanning tree of $\wh{B_G(o, n)}$ in the same number of ways. 
Letting $n$ tend to $\infty$, we see that the law of $\WSFs\cap B_G(o,m)$ is $\USF(K_m)$.
Note that $\Seq{K_m}_{m \ge 1}$ is an exhaustion of $\ol{\WSFs(G)}$. (More precisely, each component of $\ol{\WSFs(G)}$ is exhausted by the corresponding sequence  of components of  $K_m$.) 
Therefore by the definition of $\FSF$, the measures $\USF(K_m)$ converge to \(\FSF(\ol{\WSFs(G)})\) as $m\to\infty$ (restricted to any finite subgraph of $G$). 
Since the law of $\WSFs\cap B_G(o,m)$ is $\USF(K_m)$, by letting $m$ tend to $\infty$, we obtain $\WSF(G)=\FSF(\ol{\WSFs(G)})$. 
\end{proof}

\begin{proof}[Proof of Corollary~\ref{cor:wired2}]
Recall that   $\WSF(G)\preccurlyeq \FSF(G)$	for any locally finite connected graph $G$ (see, e.g., \cite[Section 10.2]{TreeBook}). 
Together with Theorem~\ref{thm:resample}, we obtain 
$$\WSF(\ol{\WSFs(G)}) \preccurlyeq \FSF(\ol{\WSFs(G)})=\WSF(G).$$
Let $(\WSFs',\WSFs)$ be a coupling of $\WSF(\ol{\WSFs(G)})$ and $\WSF(G)$ 
such that $\WSFs' \subset\WSFs$.
Since each connected component of $\WSFs'$ is an infinite graph a.s., while each component of $\WSFs$ has a single end, we must have $\WSFs=\WSFs'$ a.s.
This proves the first assertion; the second is even simpler, because by Theorem~\ref{thm:resample}, we need only show that $\WSF(G) = \FSF(G)$.
\end{proof}

\subsection{A counterexample \texorpdfstring{for $\WSF(\ol{\FSF(G)}) =\WSF(G)$}{}}\label{subsec:counter}
Recall that for any graph $G$ and neighbors $x, y$ in $G$, Kirchhoff's formula extended to the wired spanning forest gives that
$\P[(x, y) \in \WSFs(G)] = \effw^G(x, y)$; see \cite[Equation~(10.3)]{TreeBook}. Here, the
\notion{wired effective resistance}
between $x$ and $y$ is defined by
\[
\effw^G(x,y)=\Bigl(\inf\bigl\{\cE(f) \st f(x)=1,\, f(y)=0,\, \#\bigl(f^{-1}[\R \setminus \{0\}]\bigr) < \infty\bigr\}\Bigr)^{-1}.
\]
If $H$ is a subgraph of $G$ that includes $(x, y)$ but does not include at least one edge $(u, w)$
for which the wired current $\iw^{(x, y)}(u, w) \ne 0$, then Thomson's
principle (uniqueness of the minimizing-energy unit flow) yields $\P[(x, y) \in \WSFs(G)] < \P[(x, y) \in \WSFs(H)]$; see \cite[Section 9.1]{TreeBook} for the definition of wired current.

Now let $G$ be the graph consisting of two copies of $\Z^5$, which we denote by $\Z^5\times\{0\}$ and $\Z^5\times\{1\}$, and an edge $e$ connecting
$o_0\defeq (o, 0)$ and $o_1\defeq (o, 1)$. As before, $o$ represents the origin of $\Z^5$.  Since $\Z^5$ is
transient, the wired current $\iw^{e}$ is nonzero on infinitely many edges of $\Z^5 \times
\{i\}$ for each $i \in \{0, 1\}$. (In fact, it can be proved that all edges have
nonzero current.) Recall that $\FSF(\Z^5)=\WSF(\Z^5)$. Since $\FSFs(\Z^5)$
contains infinitely many trees a.s., its induced components are not all of $\Z^5$. 
Furthermore, $e$ is not contained in any cycle, whence $e \in \FSFs(G)$ a.s., and $\FSF(G)$ may be coupled with $\FSF(\Z^5 \times \{0\})$ and $\FSF(\Z^5
\times \{1\})$ so that $\FSFs(G) = \{e\} \cup \FSFs(\Z^5 \times \{0\}) \cup \FSFs(\Z^5 \times \{1\})$.
Let $\cT_e$ be the  component of $\FSFs(G)$ containing $e$. Then $\cT_e$ consists of $e$
and the component of $\FSFs(\Z^5\times \{i\})$ containing $o_i$, where $i=0,1$. Each edge of
$\Z^5$ has the same probability of being in $\ol{\FSFs(\Z^5)}$, whence infinitely many edges $(u, w)$
with $\iw^e(u, w) \ne 0$ are not in $\ol{\cT_e}$ a.s. It follows from the
preceding paragraph that 
$\P[(x, y) \in \WSFs(G)] < \P[(x, y) \in \WSFs(\ol{\cT_e}) \mid \ol{\cT_e}]$.
Taking the expectation gives the result.

This same method answers negatively a long-standing question of whether $\WSF(\FSFs(G))=\WSF(G)$
for Cayley graphs, $G$. We may take $G$ to be the natural Cayley graph of the free product of $\Z^3$
with $\Z_2$ to obtain a counterexample. The analysis is similar to the preceding.

\newcommand{\MRhref}[2]{\href{http://www.ams.org/mathscinet-getitem?mr=#1}{MR#2}}
\def\@rst #1 #2other{#1}
\newcommand\MR[1]{\relax\ifhmode\unskip\spacefactor3000 \space\fi
	\MRhref{\expandafter\@rst #1 other}{#1}}


\begin{thebibliography}{{Law}18}
	
	\bibitem[BKPS04]{GeometryofUSF}
	I.~Benjamini, H.~Kesten, Y.~Peres, and O.~Schramm.
	\newblock Geometry of the uniform spanning forest: transitions in dimensions
	{$4,8,12,\dots$}.
	\newblock {\em Ann.\ of Math.\ (2)}, 160(2):465--491, 2004. \MR{2123930
		(2005k:60026)}
	
	\bibitem[BLPS99]{BLPSgip}
	I.~Benjamini, R.~Lyons, Y.~Peres, and O.~Schramm.
	\newblock Group-invariant percolation on graphs.
	\newblock {\em Geom.\ Funct.\ Anal.}, 9(1):29--66, 1999. \MR{99m:60149}
	
	\bibitem[BLPS01]{USF}
	I.~Benjamini, R.~Lyons, Y.~Peres, and O.~Schramm.
	\newblock Uniform spanning forests.
	\newblock {\em Ann.\ Probab.}, 29(1):1--65, 2001. \MR{1825141 (2003a:60015)}
	
	\bibitem[BLS99]{BLS99}
	I.~Benjamini, R.~Lyons, and O.~Schramm.
	\newblock Percolation perturbations in potential theory and random walks.
	\newblock In M.~Picardello and W.~Woess, editors, {\em Random Walks and
		Discrete Potential Theory}, Sympos.\ Math., pages 56--84, Cambridge, 1999.
	Cambridge University Press.
	\newblock Papers from the workshop held in Cortona, 1997. \MR{1802426}
	
	\bibitem[Dub09]{Dubedat-GFF}
	J.~Dub\'{e}dat.
	\newblock S{LE} and the free field: partition functions and couplings.
	\newblock {\em J. Amer. Math. Soc.}, 22(4):995--1054, 2009. \MR{2525778}
	
	\bibitem[HN17]{Indistinguish}
	T.~Hutchcroft and A.~Nachmias.
	\newblock Indistinguishability of trees in uniform spanning forests.
	\newblock {\em Probab.\ Theory Related Fields}, 168(1-2):113--152, 2017.
	\MR{3651050}
	
	\bibitem[HP19]{Component}
	T.~Hutchcroft and Y.~Peres.
	\newblock The component graph of the uniform spanning forest: transitions in
	dimensions {$9,10,11,\ldots $}.
	\newblock {\em Probab.\ Theory Related Fields}, 175(1-2):141--208, 2019.
	\MR{4009707}
	
	\bibitem[Hut18]{Tom-End}
	T.~Hutchcroft.
	\newblock Interlacements and the wired uniform spanning forest.
	\newblock {\em Ann.\ Probab.}, 46(2):1170--1200, 2018. \MR{3773383}
	
	\bibitem[Hut19]{Hutchcroft:universality}
	T.~Hutchcroft.
	\newblock Universality of high-dimensional spanning forests and sandpiles.
	\newblock {\em Probab.\ Theory Related Fields}, Jun 2019.
	\newblock \url{http://dx.doi.org/10.1007/s00440-019-00923-3}.
	
	\bibitem[Law80]{Lawler-Bi}
	G.~F. Lawler.
	\newblock A self-avoiding random walk.
	\newblock {\em Duke Math.\ J.}, 47(3):655--693, 1980. \MR{587173}
	
	\bibitem[{Law}18]{Lawler-two-sided}
	G.~F. {Lawler}.
	\newblock {The infinite two-sided loop-erased random walk}.
	\newblock {\em ArXiv e-prints}, February 2018, 1802.06667.
	
	\bibitem[LMS08]{End}
	R.~Lyons, B.~J. Morris, and O.~Schramm.
	\newblock Ends in uniform spanning forests.
	\newblock {\em Electron. J. Probab.}, 13:paper no. 58, 1702--1725, 2008.
	\MR{2448128}
	
	\bibitem[LP16]{TreeBook}
	R.~Lyons and Y.~Peres.
	\newblock {\em Probability on trees and networks}, volume~42 of {\em Cambridge
		Series in Statistical and Probabilistic Mathematics}.
	\newblock Cambridge University Press, New York, 2016. \MR{3616205}
	
	\bibitem[LSW04]{LSW-UST}
	G.~F. Lawler, O.~Schramm, and W.~Werner.
	\newblock Conformal invariance of planar loop-erased random walks and uniform
	spanning trees.
	\newblock {\em Ann. Probab.}, 32(1B):939--995, 2004. \MR{2044671}
	
	\bibitem[LSW19]{LSW-4D}
	G.~F. {Lawler}, X.~{Sun}, and W.~{Wu}.
	\newblock {Four dimensional loop-erased random walk}.
	\newblock {\em Ann.\ Probab.}, 2019.
	\newblock To appear.
	
	\bibitem[McG91]{McGuinness:NW}
	S.~McGuinness.
	\newblock Recurrent networks and a theorem of {N}ash-{W}illiams.
	\newblock {\em J. Theoret. Probab.}, 4(1):87--100, 1991. \MR{1088394}
	
	\bibitem[Mor03]{Morris}
	B.~Morris.
	\newblock The components of the wired spanning forest are recurrent.
	\newblock {\em Probab.\ Theory Related Fields}, 125(2):259--265, 2003.
	\MR{1961344}
	
	\bibitem[NW59]{NashWil:crit}
	C.~{\relax St}. J.~A. Nash-Williams.
	\newblock Random walk and electric currents in networks.
	\newblock {\em Proc. Cambridge Philos. Soc.}, 55:181--194, 1959. \MR{23:A2239}
	
	\bibitem[Pem91]{Pemantle}
	R.~Pemantle.
	\newblock Choosing a spanning tree for the integer lattice uniformly.
	\newblock {\em Ann.\ Probab.}, 19(4):1559--1574, 1991. \MR{1127715 (92g:60014)}
	
	\bibitem[Pet83]{ergodic}
	K.~Petersen.
	\newblock {\em Ergodic theory}, volume~2 of {\em Cambridge Studies in Advanced
		Mathematics}.
	\newblock Cambridge University Press, Cambridge, 1983. \MR{833286}
	
	\bibitem[Sch00]{Schramm00}
	O.~Schramm.
	\newblock Scaling limits of loop-erased random walks and uniform spanning
	trees.
	\newblock {\em Israel J. Math.}, 118:221--288, 2000. \MR{1776084}
	
	\bibitem[She16]{Shef-burger}
	S.~Sheffield.
	\newblock Quantum gravity and inventory accumulation.
	\newblock {\em Ann. Probab.}, 44(6):3804--3848, 2016. \MR{3572324}
	
	\bibitem[Sun19]{Bernoulli-UST}
	X.~Sun.
	\newblock Random planar geometry through the lens of uniform spanning tree.
	\newblock {\em Bernoulli News}, 26(2):10--13, 2019.
	
	\bibitem[Wil96]{Wilson}
	D.~B. Wilson.
	\newblock Generating random spanning trees more quickly than the cover time.
	\newblock In {\em Proceedings of the {T}wenty-Eighth {A}nnual {ACM} {S}ymposium
		on the {T}heory of {C}omputing ({P}hiladelphia, {PA}, 1996)}, pages 296--303.
	ACM, New York, 1996. \MR{1427525}
	
\end{thebibliography}
\end{document}